\documentclass{amsart}
\usepackage[dvips]{graphicx}
\usepackage{amscd}
\usepackage{amsmath}
\usepackage{amsxtra}
\usepackage{amsfonts}
\usepackage{amssymb}
\newtheorem{theorem}{Theorem}[section]

\newtheorem{proposition}[theorem]{Proposition}
\theoremstyle{definition}
\newtheorem{definition}[theorem]{Definition}
\newtheorem{remark}[theorem]{Remark}

\newtheorem{example}[theorem]{Example}
\theoremstyle{remark}

\renewcommand{\theclaim}{\textup{\theclaim}}

\newtheorem*{acknowledgements}{Acknowledgements}

\numberwithin{equation}{section}

\newcommand{\triunghi}[7]{
\begin{table}[ht]
\begin{center}
\begin{tabular}{ccccccccc}

$ $&$ $&$ $&$ $&$#1$&$ $&$ $&$ $&$ $\\
$ $&$ $&$ $& $\nearrow $ &$ $&$\nwarrow $&$ $&$ $&$ $\\
$ $&$ $&$#4 $&$ $&$ $&$ $&$#5 $&$ $&$ $\\
$ $&$ $&$ $&$\nwarrow $&$ $&$\nearrow$&$ $&$ $&$ $\\
$ $&$\swarrow $&$ $&$ $&$#7$&$ $&$ $&$\searrow $&$ $\\
$ $&$ $&$ $&$ $&$\downarrow$&$ $&$ $&$ $&$ $\\

$#2 $&$ $&$\longleftarrow$&$ $&$#6$&$ $&$\longrightarrow$&$ $&$#3 $\\
\end{tabular}
\end{center}
\end{table}
}

\newcommand{\co}{{\emph {\bf C}}}

\newcommand{\focksum}{
\bigoplus_{{\tiny \begin{array}{c}
n\geq1\\
i_1,i_2,...,i_n\in I\\
i_1\not=i_2\not=...\not=i_n
\end{array}}
 }
 }
\newcommand{\diagrama}[6]{
\begin{table}[htdp]
\begin{center}
\begin{tabular}{cccc}
$ $&$ $&$#4$&$ $\\
$ $&$#1$&$\longrightarrow$&$#2$\\
$#5$&$\downarrow$&$\nearrow$&$ $\\
$ $&$#3$&$#6$&$ $
\end{tabular}
\end{center}
\end{table}
}

\begin{document}
\title[Generalized free amalgamated product of C*- algebras]{Generalized free amalgamated product of C*- algebras}
\author{Stefan Teodor Bildea}
\address{Department of Mathematics\\
The University of Iowa\\
14 MacLean Hall\\
Iowa City, IA 52242-1419\\
U.S.A.\\} \email{ Stefan Teodor Bildea: sbildea@math.uiowa.edu}

 \keywords{}

\begin{abstract}
We construct a generalized version for the free product of unital C*-algebras $(A_i)_{i\in I}$ with amalgamation over a family of common unital subalgebras $(B_{ij})_{i,j\in I,i\ne j}$, starting from the group-analogue. When all the subalgebras are the same, we recover the free product with amalgamation over a common subalgebra. We reduce the problem to the study of minimal amalgams. We specialize to triangles of algebras and subalgebras, study freeness in this context, and give some examples of constructions of minimal amalgams derived from triangles of operator algebras.
\end{abstract}

\maketitle \tableofcontents
\section{\label{Intro}Introduction}
The von Neumann group algebra of a free product of groups, without and with amalgamation, was the starting point in 
Dan Voiculescu's free probability theory. Due to the succes of the spatial theory of free products of C*-algebras and von Neumann algebras he developed, we find it of interest to study other kinds of amalgams, starting from what happens in the category of groups. 
In \cite{HN1}, \cite{HN2} the generalized free product with amalgamated subgroups is introduced and studied extensively. Given a family of groups $(G_i)_{i\in I}$ and subgroups $H_{ij}\subset G_i,\forall j,i\in I,j\not= i $ so that $H_{ij}\stackrel{\varphi_{ij}}{\simeq} H_{ji},\forall i,j\in I, i\not=j$, and $\varphi_{ji}=\varphi_{ij}^{-1}$, the generalized free product G of the $G_i's$ with amalgamated subgroups $H_{ij}$ is defined as follows. Take the free product of the $G_i's$ and factor it by the normal subgroup $F$ generated by elements of the form 
\[ \varphi_{ij}(h_{ij})h_{ji}^{-1}, \quad h_{ij}\in H_{ij},h_{ji}\in H_{ji}\quad i,j\in I,i\not= j,\] 
i.e. the elements of the generalized free product of groups are words over the "alphabet" $\cup_{i\in I}G_i$such that $.....(g_ih_{ij})g_j...=...g_i(\varphi_{ij}(h_{ji})g_j)...$ .
Unlike for the free product of groups, where the groups embed in their free product, it is not necessary the case that the groups $G_i$ have isomorphic copies inside of  $G$. When this is true we say that the generalized free product is realizable. 
Examples of collapsing families of groups can be found in \cite{HN1}(Example 3.2), \cite{BHN}(Chapter II,11) and more recently in \cite{Tits}. One striking sufficient condition for the existence of a realizable generalized free amalgamated product of a triangle of groups as the one bellow was given by John Stallings in \cite{Sta}. 
\triunghi{G_1}{G_2}{G_3}{H_{12}}{H_{23}}{H_{13}}{K}
He defined an angle between the subgroups $H_{12},H_{13}\subset G_1$ as follows. Consider the unique morphism $\phi:H_{12}\ast_K H_{13}\to G_1$. Let $2n=\inf \{|w| | w\in ker\phi,w\ne e\} $ and define $\theta_1=\frac{\pi}{n}$ to be the angle between the two subgroups. Consider the 3 angles involved in the triangle above. If $\theta_1+\theta_2+\theta_3\le \pi$ then the considered generalized free product is realizable. On the other hand, a famous example of a group that can be presented as the minimal amalgam of a triangle of finite groups and which does not satisfy this condition is Thompsons's infinite simple group $G_{2,1}$ (see \cite{Brown}).

\section{\label{free prod groups}Generalized free product of groups with amalgamated subgroups}
\par

\begin{definition}
Let $(G_i)_{i\in I}$ be a family of groups and for each $i\in I$ let $(H_{ij})_{j\in I,j\ne i }$ be a family of subgroups of $G_i$ as above. The generalized group free product of the family $(G_i)_{i\in I}$ with amalgamated subgroups $(H_{ij})_{i,j\in I,i\ne j}$ is the unique (up to isomorphism) group $G$ together with homomorphisms $\psi_i:G_i\to G$ satisfying the properties:
\newline (i) the diagram bellow commutes for every $i,j\in I,i\not=j;$
\diagrama{H_{ij}}{G}{H_{ji}}{\psi_i}{\varphi_{ij}}{\psi_j}
\newline (ii) given any group $K$ and homomorphisms $\phi_i:G_i\to K$ with commutative diagrams
\diagrama{H_{ij}}{K}{H_{ji}}{\phi_i}{\varphi_{ij}}{\phi_j}
\newline for every $i,j\in I,i\not=j$, there exists a unique homomorphism $\Phi:G\to K$ so that the following diagram commutes:
\diagrama{G_i}{G}{K}{\psi_i}{\phi_i}{\Phi}
\newline
We will use the notation \[ G=\ast_{i\in I}(G_i,(H_{ij})_{j\ne i}) .\]
When $G_i=\bigvee_{j\in I,j\ne i}H_{ij}$, we call $G$ the minimal generalized free product or the minimal amalgam of the family $(G_i)_{i\in I}$.
\end{definition}

\begin{remark}
We will also use the following alternative description of the family of groups and subgroups that we will amalgamate : let $G_i$ be groups for $i\in I$ and for each $i,j\in I,i\ne j$ let $H_{ij}(=H_{ji})$ be groups and $\varphi_{ij}:H_{ij}\to G_i$ be injective group homomorphisms. The groups $H_{ij}$ are now abstract groups that have isomorphic copies inside the bigger groups $G_i$ and $G_j$, and this copies will be identified. If we let 
$ K_{ij}:=\varphi_{ij}(H_{ij}) \stackrel{\varphi_{ji}\circ \varphi_{ij}^{-1}}{\simeq} K_{ji}:=\varphi_{ji}(H_{ij})$,
then the family of groups and subgroups $(G_i)_{i\in I}$ , $(K_{ij})_{j\in I,j\ne i}$ is as in the previous definition. By the generalized free product of the family $(G_i)_{i\in I}$ with amalgamated subgroups $(H_{ij})_{j\in I ,j\ne i}$ we will mean $\ast_{i\in I}(G_i,(K_{ij})_{j\in I,j\ne i})$.  Of course we can go both ways with the description of the families and we will use both notation alternatively in this paper.
\end{remark}

\par The first criterion for the existence of the isomorphic copies of the groups $G_i$ inside $G$ refers to the reduced amalgam. For $i\in I$ let $H_i$ be the subgroup of $G_i$ generated by the family $(H_{ij})_{j\in I,j\ne i}$, i.e. $H_i:=\bigvee_{j\in I,j\ne i}H_{ij} \subseteq  G_i$. The following reduction theorem is due to Hanna Neumann (\cite{HN1}) and has also a nice presentation in \cite{BHN}, section 15. We present the proof for consistency.

\begin{theorem}
With the above notations, the generalized free product G of the family $(G_i)_{i\in I}$ with amalgamated subgroups $(H_{ij})_{i,j\in I,j\ne i}$ is realizable if and only if the generalized free product H of the family $(H_i)_{i\in I}$ with amalgamated subgroups $(H_{ij})_{i,j\in I,j\ne i})$ is realizable.
\end{theorem} 

{\it Proof:   }
It is clear that if $G$ is realizable, then each $H_i$ will have isomorphic copies as subgroups of the isomorphic copies of $G_i$, hence $H$ is obviously realizable. Suppose now that $H$ is realizable. For each $i\in I$ we first take the free product  $G_i\ast_{H_i}H$.  Next consider 
\[ G=\ast_{H,i\in I}(G_i\ast_{H_i}H).\]
Then each $G_i$ has an isomorphic copy in $G$. Furthermore, since $H$ is the amalgamated subgroup,  and it is the realizable minimal amalgam, we have, using the same notation for the copies of $G_i$ inside $G$ :  $G_i\cap G_j=G_i\cap G_j\cap H=H_i\cap H_j=H_{ij}=H_{ji}\subset G$. It follows that the $G_i's$ embed in $\bigvee_{i\in I}\subset G$ such that the resulting diagrams concerning the subgroups commute, hence $\ast_{i\in I}(G_i,(H_{ij})_{j\ne i})$ is realizable. 
\newline
Notice that actually $G=\bigvee_{i\in I}G_i$. Indeed, a moments thought suffices to see that an element $x\in G$ has the form
\[ x=w_1(\{H_i|i\in I\})g_1w_2(\{H_i|i\in I\})g_2...w_n(\{H_i|i\in I\})g_n,\]
with $n\geq 1$ and where $g_k\in G_{i_k},i_1\ne i_2\ne...\ne i_n$ and $w_k(\{H_i | i\in I\})$ are words over the "alphabet" $\cup_{i\in I}H_i$ so that the letters adjacent to an element of $G_{i_k}$ do not belong to $\cup_{j\in I,j\ne i_k}H_{i_k,j}$. But $w_k(\{H_i | i\in I\})\in \bigvee_{i\in I}G_i$, hence $x\in \bigvee_{i\in I}G_i$.
\rule{5 pt}{5 pt}

\begin{remark}
We will be interested in this paper in a particular kind of family of groups, that can be represented on a triangle diagram.
\triunghi{G_1}{G_2}{G_3}{H_{12}}{H_{13}}{H_{23}}{H_{123}}
All arrows are injective group homomorphisms. The addition of the group $H_{123}$ is just a necessary condition for the generalized free product to be realizable. Indeed, if $\ast_{i\in I}(G_i,(H_{ij})_{i\ne j})$, $I=\{ 1,2,3\}$, is realizable, then we may assume that the arrows are inclusions and that $G_i\cap G_j=H_{ij}$. Then $G_1\cap G_2\cap G_3=H_{12}\cap H_{13}=H_{12}\cap H_{23}=H_{13}\cap H_{23}=H_{123}$. Therefore we will assume the triangles to be fillable, i.e. the intersections of pairs of edge groups are isomorphic under the given family of maps; for example 
\[ H_{123}\simeq \varphi_{12}(H_{12})\cap \varphi_{13}(H_{13})(\subset G_1)\stackrel{\varphi_{31}\circ \varphi_{13}^{-1}}{\simeq}\varphi_{32}(H_{23})\cap \varphi_{31}(H_{13})\subset G_3.\]
\end{remark}

\begin{definition} Suppose we are given a triangle of groups and injective group morphisms as the one above, with the additional hypothesis that the vertex groups are generated by the pairs of images of the edge groups - we will call such a triangle a minimal triangle. If the corresponding generalized free product with amalgamations exists, we say the triangle is realizable.
\end{definition}

The last result of this section addresses the realizability of triangles of groups as the one above. 
\begin{theorem}\label{retri}
Given a minimal triangle of groups and injective groups homomorphisms as the previous one, each of the following is a sufficient condition for the triangle to be realizable:
\newline
(i) One of the vertex groups, $G_3$ say, is the free product of $H_{23},H_{13}$ with amalgamation over $H_{123}$, and with inclusion mappings $\varphi_{32},\varphi_{13}$. In this case we have $\ast_{i=1}^3(G_i,(H_{ij})_{j\ne i})=G_1\ast_{H_{12}}G_2$. 
\newline
(ii) Two of the the groups $G_i$, $G_1$ and $G_2$ say, have the property that every element of $H_{12}$ commutes with every element of $H_{13}$ and of $H_{23}\quad (i=1,2)$. In this case $\ast_{i=1}^{3}(G_i,(H_{ij})_{j\ne i})$ is a quotient of $G_1\ast_{H_{12}}G_2$. 
\end{theorem}

This theorem and its proof can be found in \cite{HN1}, section 9. The idea of the proof is to look at the way $H_{13}\ast_{H_{123}}H_{23}$ embeds into $G_1\ast_{H_{12}}G_2$. This is the way to prove similar results for triangles of operator algebras.

\section{\label{s0}Full and reduced products of unital C*-algebras}
This section contains a brief review on full and reduced free products of unital C*-algebras with amalgamation over a common subalgebra.

\begin{definition}
Given a family of unital C*-algebras $(A_i)_{i\in I}$ ($I$ is a set having at least two elements) with a common unital C*-subalgebra $B$, and injective $\ast$-homomorphisms $\phi_i:B\to A_i,\forall i\in I$, the corresponding full amalgamated free product C*-algebra is a unital C*-algebra $A$, equipped with injective $\ast$-homomorphisms $\sigma_i:A_i\to A,\forall i\in I$ such that :
\newline (i) $\sigma_i\circ \phi_i=\sigma_j\circ \phi_j$ for all $i,j\in I$ and $A$ is generated by $\cup \sigma_i(A_i)$;
\newline (ii) for any unital C*-algebra $D$ and any injective $\ast$-homomorphisms $\pi_i:A_i\to D$, $i\in I,$ satisfying $\pi_i\circ \phi_i=\pi_j\circ \phi_j$ for all $i,j\in I$, there is a $\ast$-homomorphism $\pi:A\to D$ such that the next diagram commutes.
\diagrama{A_i}{D}{A}{\pi_i}{\phi_i}{\pi}
\end{definition}

\begin{remark}
If $B$ is just the field of complex numbers, then the injectivity of the maps $\sigma_i$ can be omitted from the definition and follows from the fact that the full free product has enough representations ( see \cite{Picu} ). For the general case, injectivity cannot be omitted from the definition without changing the category of objects. If, for example, for each $i\in I$ there exists a conditional expectations $E_i:A_i\to B$ with faithful GNS representation (see \ref{fGNS}), then injectivity follows using the properties of the reduced free product with amalgamation and the universal property of the full free product. 
\end{remark}
\par In order to give criteria for existence of generalized amalgamated free products, we need  the following results on embeddings of full and reduced amalgamated free product C*-algebras. This first proposition works actually also for free products of non-unital C*-algebras (see \cite{ExelDy} proposition 2.4).

\begin{proposition}\label{exel} Suppose\newline 

\begin{table}[htbp]
\begin{center}
\begin{tabular}{ cccccc}

$\tilde{A} $&$\hookleftarrow $ &$\tilde{D}$&$\hookrightarrow$&$\tilde{B}$\\
$\uparrow  $&$$&$\uparrow$&$$&$\uparrow $\\
$A$&$\hookleftarrow $ &$D$&$\hookrightarrow$&$B$
\end{tabular}
\end{center}
\end{table}
is a commuting diagram of inclusions of C*-algebras. Let $\lambda:A\ast_DB\to \tilde{A}\ast_{\tilde{D}}\tilde{B}$ be the resulting $\ast$-homomorphism of full free product C*-algebras. Suppose there are conditional expectations $E_A:\tilde{A}\to A,E_D:\tilde{D}\to D$ and $E_B:\tilde{B}\to B$ onto $A,D$ and $B$, respectively, such that the diagram

\begin{table}[htbp]
\begin{center}
\begin{tabular}{ cccccc}

$\tilde{A} \quad$&$\hookleftarrow $ &$\tilde{D}\quad$&$\hookrightarrow$&$\tilde{B}\quad$\\
$\downarrow \mbox{{\sc e}}_A $&$$&$\downarrow \mbox{{\sc e}}_D$&$$&$\downarrow \mbox{{\sc e}}_B$\\
$A\quad$&$\hookleftarrow $ &$D\quad$&$\hookrightarrow$&$B\quad$
\end{tabular}
\end{center}
\end{table}

commutes. Then $\lambda$ is injective.

\end{proposition}

\par  Suppose $\phi :A \to B$ is a conditional expectation and let $E =L^2(A ,\phi )$ be the (right) Hilbert B-module obtained from $A$ by separation and completion with respect to the norm $\| a\|=\| (a,a)_{E }\|^{\frac12}$, where $(\cdot, \cdot)_{E }$ is the $B$-valued inner product, $(a_1, a_2)_{E }=\phi (a_1^{\ast}a_2)$. We denote the map $A \to E $ arising from from the definition by $a\mapsto \hat{a}$. 

Let $\pi :A \to \mathcal{L}(E )$ denote the $\ast$-representation defined by $\pi (a)\hat{b}=\hat{ab}$, where $\mathcal{L}(E )$ is the C*-algebra of all adjointable bounded $B$-module operators on $E $. Consider the specified element $\xi=\widehat{1_A}\in E$. 

\begin{definition}\label{fGNS} With the above notations, we will call $(\pi, E, \xi)$ the GNS representation of $(A,\phi)$ and write this as
\[ (\pi, E, \xi)=GNS(A, \phi).\]  
We will say that the conditional expectation $\phi:A\to B$ has faithful GNS representation if $\pi$ is a faithful representation.
\end{definition}

\begin{remark}
Note that the faithfulness of $\pi$ is equivalent to the following condition:
\[ \forall a\in A\setminus \{0\} \quad \exists x\in A\quad  \phi(x^{*}a^{*}ax)\ne0.\]
\end{remark}

\par Let $B$ be a unital C*-algebra, let $I$ be a set having at least two elements and for every $i\in I$ let $A_i$ be a unital C*-algebra containing a copy of B as a unital C*-subalgebra.
Voiculescu's construction of the reduced amalgamated free product $(A,\phi)$ starts with the Hilbert $B$-modules $(\pi_i, E_i, \xi_i)=GNS(A_i, \phi_i)$.  Next he constructs a Hilbert $B$-module $E$ which turns out to conincide with $L^2(A, \phi)$. Letting $E^{\circ}_i=L^2(A_i, \phi_i)\ominus \xi B$, define

\[ E=\xi B\oplus \focksum E^{\circ}_{i_1}\otimes_B E^{\circ}_{i_2}\otimes_B...\otimes_B E^{\circ}_{i_n}.\]
Here $\xi =1_B$ and $\xi B$ is the C*-alegbra $B$ viewed as a Hilbert $B$-module with specified element $\xi$ and all the tensor products are internal, arising from $\ast$-homomorphisms $P^{\circ}_i\pi_i|_BP^{\circ}_i$ from $B$ to $\mathcal{L}(E^{\circ}_i)$ ($P^{\circ}_i$ is the projection onto $E^{\circ}_i$). The Hilbert module $E$ is called the free product of the $E_i$ with respect to the specified vectors $\xi_i$, and is denoted by $(E,\xi)=\ast_{i\in I}(E_i, \xi_i)$.
Next, we define the free product of the $\ast$-representations $\pi_i$. For $i\in I$ let
\[ E(i)=\eta_iB\oplus 
\bigoplus_{{\tiny \begin{array}{c}
n\geq1\\
i_1,i_2,...,i_n\in I\\
i_1\not=i_2\not=...\not=i_n\\
i_1\ne i
\end{array}}}
E^{\circ}_{i_1}\otimes_B E^{\circ}_{i_2}\otimes_B...\otimes_B E^{\circ}_{i_n},\]
with $\eta_i=1_B$ and $\eta_iB$ a copy of the Hilbert $B$-module $B$, and let \[ V_i:E_i\otimes_BE(i)\to E \] be the unitary defined as follows :

\[ \begin{array}{c} \xi_i\otimes \eta_i \stackrel{V_i}{\mapsto} \xi;\\
\zeta \otimes \eta_i \stackrel{V_i}{\mapsto} \zeta; \\
\xi_i\otimes (\zeta_1 \otimes\zeta_2\otimes  ...\zeta_n)\stackrel{V_i}{\mapsto} \zeta_1 \otimes\zeta_2\otimes  ...\zeta_n;\\
\zeta \otimes (\zeta_1 \otimes\zeta_2\otimes  ...\zeta_n)\stackrel{V_i}{\mapsto}\zeta\otimes  \zeta_1 \otimes \zeta_2\otimes  ...\zeta_n,\\
\end{array}
\]
$\forall \zeta \in E^{\circ}_i, \forall \zeta_j\in E^{\circ}_{i_j},\quad i \ne i_1\ne i_2\ne...\ne i_{n-1}\ne i_n$.
Let $\lambda_i:A_i\to \mathcal{L}(E)$ be the $\ast$-homomorphism given by
\[ \lambda_i(a)=V_i(\pi_i(a)\otimes 1)V_i^{\ast}.\]
Then $A$ is defined to be the C*-algebra generated by $\cup_{i\in I}\lambda_i(A_i)$, and $\phi:A\to B$ is the conditional expectation $\phi(\cdot)=(\xi, \cdot \xi)_E$. The pair  $(A,\phi)$, together with the embeddings $\lambda_i:A_i\hookrightarrow A$ which restrict to the identity on $B$, is called the reduced amalgamated free product of the $(A_i,\phi_i)'s$ and is characterized by the following properties:
\newline (i) $\forall i\in I\quad \phi|_{A_i}=\phi_i$;
\newline (ii) the family $(A_i)_{i\in I}$ is free with respect to $\phi$, i. e.
\[ \phi(a_1a_2...a_n)=0 \quad \forall a_j\in A_{i_j}\cap\mbox{ker}\phi_{i_j},\mbox{ with }i_1\ne i_2\ne...\ne i_n,n\geq 1.\]
\newline (iii) $A$ is generated by $\cup_{i\in I}\lambda(A_i)$;
\newline (iv) the GNS representation of $\phi$ is faithful on $A$. We will write
\[(A,\phi)=\ast_{i\in I}(A_i,\phi_i).\]

As in the case of the tensor product, for the reduced free product we can construct canonical conditional expectations into the factors.
The following proposition can be found in \cite{DyBl} as Lemma 1.1:

\begin{proposition}\label{bl}
 Let $( A_i)_{i\in I}$ be a family of unital C*-algebras having a common subalgebra $1\in B\subset A_i,\forall i\in I$. Suppose there are conditional expectations $\phi_i:A\to B$, for all $i\in I$ whith faithful GNS representation. If 
 \[(A, \phi)=\ast_{i\in I}(A_i, \phi_i)\] is the reduced amalgamated free product, then for every $i_0\in I$ there is a canonical conditional expectation, $\psi_{i_0}:A\to A_{i_0}$ with the following properties:

\begin{itemize}

 \item[(i)] $\psi_{i_0}|_{A_i}=\phi_i$ for every $i\in I\setminus \{ i_0 \}$ \\
 \item[(ii)] $\psi_{i_0}(a_1a_2...a_n)=0$ whenever $n\geq 2$ and $a_j\in A_{i_j}\cap ker\phi_{i_j}$ with $i_1\not=i_2\not=...\not=i_{n-1}\not=i_n$.

\end{itemize}
\end{proposition}

This last theorem we mention refers to the embedding of reduced free products of unital C*-algebras and can be found also in \cite{DyBl} as theorem 1.3.
\begin{theorem}\label{redfreeprod}
Let $B\subset \tilde{B}$ be a (not necessarily unital) inclusion of unital C*-algebras. Let $I$ be a set and for each $i\in I$ suppose\\

\begin{table}[htbp]
\begin{center}
\begin{tabular}{cccc}

$1_{\tilde{A_i}}\in $&$\tilde{B}$&$\subset$&$ \tilde{A_i}$\\
$$&$\cup$&$ $&$\cup$\\
$1_{A_i}\in $&$B$&$\subset$&$A$
\end{tabular}
\end{center}
\end{table}
are inclusions of C*-algebras. Suppose that $\tilde{\phi_i}:A_i\to B_i$ is a conditional expectation such that $\tilde{\phi_i}(A_i)\subset B$ and assume that  $\tilde{\phi_i}$ and the restriction $\tilde{\phi_i}|_{A_i}$ have faithful GNS representations, for all $i\in I$. Let
\[(\tilde{A}, \tilde{\phi})=\ast_{i\in I} (\tilde{A_i} \tilde{\phi_i})\]
\[ (A,\phi)=\ast_{i\in I}(A_i,\tilde{\phi_i}|_{A_i})\]
be the reduced amalgamated free products of C*-algebras. Then there is a unique $\ast$-homomorphism ${\it k}:A\to \tilde{A}$ such that for every $i\in I$ the diagram\\

\begin{table}[htbp]
\begin{center}
\begin{tabular}{ccc}

${\tilde{A_i}} $&$\hookrightarrow$&$ \tilde{A_i}$\\
$\cup$&$ $&$\uparrow {\it k}$\\
$A_i$&$\hookrightarrow$&$A$
\end{tabular}
\end{center}
\end{table}
commutes, where the horizontal arrows are the inclusions arising from the free product construction. Moreover, ${\it k}$ is necessarily injective.

\end{theorem}

\section{\label{s1}Generalized free products of algebras with amalgamated subalgebras}
\par

\begin{definition} Let $(A_i))_{i\in I}$ be a family of unital algebras, for all $i\in I$ let $B_{ij} \subset A_i$ for all $j\in I,j\ne i$ be unital subalgebras such that $B_{ij}\stackrel{\phi_{ij}}{\simeq}B_{ji}$, $\phi_{ij}^{-1}=\phi_{ji},j\in I,j\ne i$. The generalized free product of the family $(A_i)_{i\in I}$ with amalgamated subalgebras $(B_{ij})_{i,j\in I,i\ne j}$ is the unique unital algebra $A$ (up to isomorphism) together with unital homomorphisms $\psi_i:A_i\to A$ satisfying the properties:
\begin{itemize}
\item[(i)] the diagram 
\diagrama{B_{ij}}{A}{B_{ji}}{\psi_i}{\phi_{ij}}{\psi_j}
\newline commutes for all $i,j\in I,i\ne j$ and $A$ is generated by $\cup_{i\in I}\psi_i(A_i)$.\\
\item[(ii)] for any unital algebra $C$ and unital homomorphims $\varphi_i:A_i\to C$ that make the diagram
\diagrama{B_{ij}}{C}{B_{ji}}{\phi_i}{\phi_{ij}}{\phi_j}
\\ commute, there exists a unique homomorphism $\Phi:A\to B$ making the diagram bellow commute.
\diagrama{A_{i}}{A}{C}{\psi_i}{\phi_{i}}{\Phi}
\end{itemize}
We will use the notation:
\[ A=\ast_{i\in I}(A_i,(B_{ij})_{j\ne i}) .\]
When $A_i=\bigvee_{j\in I,j\ne i}B_{ij}$, we call $A$ the minimal generalized free product, or the minimal amalgam.
\end{definition}
\begin{remark}
Similarly to case of groups, we do not expect the $\ast$-homomorphisms $\psi_i:A_i\to A$ to be always injective - passing from collapsing family of groups to their group algebras, we get collapsing families of algebras. When the $\psi_i^{,}$-s are injective, we say that the generalized free product is realizable.
\end{remark}
\begin{example}\label{dmodul}
Suppose $I=\{1,2,3\}$ and $\ast_{i=1}^3(A_i, (B_{ij})_{j\ne i })$ is realizable. We can present the family of unital algebras $(A_i)_{i\in I}$, with unital subalgebras $(B_{ij})_{ij}$ as in the previous definition, in the form of a triangle of algebras. First, let us remark that, if $A=\ast_{i\in I}(A_i,(B_{ij})_{j\in I,j\ne i})$ is the {\bf realizable} generalized free product, with $\psi_i:A_i\to A$ the corresponding injective homomorhpisms, then \[ \cap_i(\psi_i(A_i))=\psi_1(B_{12})\cap \psi_2(B_{23})=\psi_2(B_{23})\cap \psi_3(B_{13}) ,\]
hence letting $D=\psi_1^{-1}(\cap_{i\in I}(\psi_i(A_i)))$, we have $D\subset B_{ij}$, for all $i,j\in I,i\ne j$. 
We get a  triangle of algebras and inclusions, which will be denoted \\$(A_1,A_2,A_3;B_{12},B_{13},B_{23};D)$ - we will employ this notation also for triangles of groups or C*-algebras.
\triunghi{A_1}{A_2}{A_3}{B_{12}}{B_{13}}{B_{23}}{D}
Therefore, for such families, a necessary condition for realizability is that the intersections of the edge algebras in the vertex algebras have to be isomorphic and identified in the final amalgam. Such triangles are called fillable. We will consider only fillable triangles.
As a $D$ module, $A$ is the quotient of the $D$-module which has as basis the set
\[ B=\{a_1a_2...a_n | n\in {\emph {\bf N}},a_j\in A_{i_j},i_1\ne i_2\ne ...\ne i_n\} \]
by the submodule generated by relations of the form
\[ a_1...a_{j-1}\left(\lambda a_j^{(0)}+\mu a_j^{(1)}\right)a_{j+i}...a_n \]
\[=\lambda a_1...a_{j-1}a_j^{(0)}a_{j+1}...a_n+\mu a_1...a_{j-1}a_j^{(1)}a_{j+1}...a_n\mbox{ where }\lambda ,\mu \in \co,\]
\[a_j=1\Rightarrow a_1...a_n=a_1...a_{j-i}a_{j+i}...a_n,\]
\[ a_{j-1}\in A_k,a_j\in B_{kl},a_{j+1}\in A_l\Rightarrow a_1...(a_{j-1}a_j)a_{j+1}...a_n=a_1...a_{j-1}(\phi_{kl}(a_j)a_{j+1})...a_n\]
 where $k,l\in I,k\ne l$ and $A_k\supseteq B_{kl}\stackrel{\phi_{kl}}{\simeq}B_{lk}\subseteq A_l$.
\end{example}

\begin{remark}
If each of the $A_i$ has an involution, then so does $A$ in the obvious manner and, moreover, $A$ is the universal object as in the previous definition in the category of unital $\ast$-algebras and unital $\ast$-homomorphisms.
\end{remark}

The following proposition relates the generalized free products of groups and algebras. If $G$ is a group, then its group algebra $\co[G]$ is the $\ast$-algebra with basis $G$, multiplication given by group multiplication and involution by $g$*$=g^{-1}$.
\begin{proposition}\label{groupalg}(Group Algebras)
Let $G$ be the generalized free product of the groups $(G_i)_{i\in I}$ with amalgamated subgroups $(H_{ij})_{i,j\in I,j\ne i}$. Then \co$[G]$ is the generalized free product of $(\co[G_i])_{i\in I}$ with amalgamation over $(\co[H_{ij}])_{i,j\in I,j\ne i}$.
\end{proposition}

{\it Proof.} By the definition of the generalized free product of groups,  there are group morphims $f_i:G_i\to G$; this extend by linearity to $\ast$-algebra morphisms $f_i:\co[G_i] \to \co[G]$. To show that $\co[G]$ is the universal object described above, let $B$ be a fixed unital $\ast$-algebra and let $\phi_i:\co[G_i]\to B$ be unital $\ast$-morphisms such that the diagrams from the definition, with $B_i=\co[G_i]$, etc., commute. A morfism $\phi:\co[G]\to B$ is completely determined by its restriction $\phi|_0:G\to \mathcal{U}(B)$, where $\mathcal{U}(B)$ are the unitaries of $B$. Using the definition of $G$ and the restrictions $\phi_i|_{G_i}:G_i\to \mathcal{U}
(B)$, we obtain the desired restriction $\phi_0:G\to \mathcal{U}(B)$.
\rule{5 pt}{5 pt}

The generalized free product for unital C*-algebras and the corresponding minimal amalgam is the universal object defined similarly to that for unital algebras. 

\begin{definition} Let $(A_i))_{i\in I}$ be a family of unital C*-algebras, for all $i\in I$ let $B_{ij} \subset A_i$ for all $j\in I,j\ne i$ be unital C*-subalgebras such that $B_{ij}\stackrel{\phi_{ij}}{\simeq}B_{ji}$, $\phi_{ij}^{-1}=\phi_{ji},j\in I,j\ne i$. The generalized free product of the family $(A_i)_{i\in I}$ with amalgamated subalgebras $(B_{ij})_{i,j\in I,i\ne j}$ is the unique unital C*-algebra $A$ (up to C*-isomorphism) together with unital $\ast$- homomorphisms $\psi_i:A_i\to A$ satisfying the properties:
\begin{itemize}
\item[(i)] the following diagram 
\diagrama{B_{ij}}{A}{B_{ji}}{\psi_i}{\phi_{ij}}{\psi_j}
commutes for all $i,j\in I,i\ne j$ and $A$ is generated by $\cup_{i\in I}\psi_i(A_i)$.\\
\item[(ii)] for any unital C*-algebra $D$ and unital $\ast$-homomorphisms $\varphi_i:A_i\to D$ that make the diagram
\diagrama{B_{ij}}{D}{B_{ji}}{\varphi_i}{\phi_{ij}}{\varphi_j}
commute, there exists a unique $\ast$-homomorphism $\Phi:A\to D$ making the diagram bellow commute.
\diagrama{A_{i}}{A}{D}{\psi_i}{\varphi_i}{\Phi.}
\end{itemize}
We will use the notation 
\[ A=\ast_{i\in I}(A_i,(B_{ij})_{j\ne i}) .\]
When $A_i=\bigvee_{j\in I,j\ne i}B_{ij}$, we call $A$ the minimal generalized free product, or the minimal amalgam.
\end{definition}

\begin{remark}
It is easy to see that such a C*-algebra exists - just take the enveloping C*-algebra of the $\ast$-generalized free product of the family $(A_i)_{i\in I}$. Unlike the free product, the $\ast$-homomorphisms $\psi_i:A_i\to A$ need not be injective - such an example relies on the examples with collapsing families of groups and the next proposition which gives the link between groups and C*-algebras of groups.
\end{remark}

\begin{remark}
In what follows, all groups will be discrete. If $G$ is a (discrete) group, we will denote by $C^{*}(G)$ the group C*-algebra of $G$, and with $C^{*}_{red}(G)$ the reduced group C*-algebra of $G$.
\end{remark}
\begin{proposition}\label{gpcstar}(Group C*-algebras)
Let $G$ be the generalized free product of the groups $(G_i)_{i\in I}$ with amalgamated subgroups $(H_{ij})_{i,j\in I,j\ne i}$. Then the group C*-algebra  C*$(G)$ is the generalized free product of C*$(G_i)_{i\in I}$ with amalgamation over C*$(H_{ij})_{i,j\in I,j\ne i}$.
\end{proposition}
{\it Proof.} Using the above remark, the definition of the generalized free product, and, as in the previous proof,  the fact that a unitary representation of the generalized free product of groups is determined by its restriction to the factors $G_i$ , the proposition follows immediately.  
\rule{5 pt}{5 pt}

\begin{remark}\label{realize}
To show that a generalized free product $A$ of a family $(A_i)_{i\in I}$ of  unital algebras ( $\ast$-algebras, C*-algebras ) with amalgamated subalgebras $(B_{ij})_{j\in I,j\ne i}$ is realizable, the injectivity of the homomorphisms $\psi:A_i\to A$ has to be proved. This will follow via the universal property  once we find an algebra ($\ast$-algebra, C*- algebra) $B$ and injective ($\ast$-) homomorphisms  $\phi_i:A_i\to B$ so that the following diagram commutes for each $i,j\in I,i\ne j$ :
\newline
\diagrama{B_{ij}}{B}{B_{ji}}{\phi_i}{\phi_{ij}}{\phi_j.}
\end{remark}


We address now the question of realizability for generalized free products of algebras.  For group C*-algebras this amounts to the realizability of the generalized free product of groups.

\begin{theorem}
Let $(G_i)_{i\in I}$ be groups and for each $i\in I$ and for each $j\in I,j\not=i $ let $H_{ij}$ be a subgroup of $G_i$ such that $H_{ij}\stackrel{\phi_{ij}}{\simeq} H_{ji}$, with $\phi_{ij}^{-1}=\phi_{ji},\forall i,j\in I,i\not=j$. If the generalized free product of $(G_i)_{i\in I}$ with amalgamated subgroups $(H_{ij})_{i,j\in I,i\not=j}$ is realizable then :\\
(i) the generalized free product of the family
of group algebras $(\co[G_i])_{i\in I}$ with amalgamated subalgebras  $(\co[H_{ij}])_{i,j\in I,i\not=j}$ is realizable;\\
(ii) the generalized free product of the family
of group C*-algebras $C^{*}(G_i)_{i\in I}$ with amalgamated subalgebras  $C^{*}(H_{ij})_{i,j\in I,i\not=j}$ is realizable.
\end{theorem}

{\it Proof}
(i) From propostion \ref{groupalg} we know that $\co[G]=\ast_{i\in I}(\co[G_i],(\co[H_{ij}])_{j\ne i})$. The corresponding maps $f_i:\co[G_i]\to \co[G]$ are the extensions of the injective group homomorphism $\psi_i:G_i\to G$ from the definition of $\ast_{i\in I}(G_i,(H_{ij})_{j\ne i})$ and therefore they are also injective and realizability is proved. \\
(ii) Use \ref{gpcstar} and the fact that the inclusion $G_i\subset G$ induces an inclusion $C^{*}(G_i)\subset C^{*}(G)$.
\rule{5 pt}{5 pt}


As in the case of groups, triangles of unital C*- algebras are of special interest. They are the simplest non-trivial case of a generalized free product with amalgamation. Next we present a version of theorem \ref{retri} for C*-algebras. 
\begin{theorem}
Let $(A_1,A_2,A_3;B_{12},B_{13},B_{23};B_{123})$ be a minimal triangle of unital C*-algebras (i.e. the vertex algebras are generated by the images of the edge algebras).  Each of the following is a sufficient condition for the realizability of the the generalized free product with amalgamation $\ast_{i=1}^3(A_i,(B_{ij})_{j\ne i})$:
\newline (i) $A_1$ is the full free product of $B_{12}$ and $B_{13}$ with amalgamation over $B_{123}$ and there exist conditional expectations $E_{12}:A_2\to B_{12},E_{13}:A_3\to B_{13}$ and $E_{123}:B_{23}\to B_{123}$ such that the diagram
\newline 
\begin{table}[htbp]
\begin{center}
\begin{tabular}{ cccccc}

$A_2\quad $&$\hookleftarrow $ &$B_{23}\quad$&$\hookrightarrow$&$A_3\quad$\\
$\downarrow \mbox{{\sc e}}_{12}  $&$$&$\downarrow \mbox{{\sc e}}_{123}$&$$&$\downarrow \mbox{{\sc e}}_{13}$\\
$B_{12}\quad $&$\hookleftarrow $ &$B_{123}\quad$&$\hookrightarrow$&$B_{13}\quad$
\end{tabular}
\end{center}
\end{table}

commutes. 
\newline (ii) There exist conditional expectations $E^{B_{1i}}_{B_{123}}:B_{1i}\to B_{123},i=2,3$, with faithful GNS representations, $A_3$ is the reduced free product of $(B_{12},E^{B_{12}}_{B_{123}})$ and $(B_{13},E^{B_{13}}_{B_{123}})$ with amalgamation over $B_{123}$ and there are conditional expectations $E^{A_i}_{B_{23}} :A_i\to B_{23}$ such that $E^{A_i}_{B_{23}}(B_{1i})\subseteq B_{123}$ and both $E^{A_{i}}_{B_{23}}$ and its restriction to $B_{1i}$ have faithful GNS representations ($i=2,3$).
\end{theorem}

{\bf Proof: } We show that the vertex algebras embed respectively into the full free product $A_2\ast_{B_{23}}A_3$ for $(i)$ and the reduced free product $(A_2\ast_{B_{23}}A_3,E^{A_2}_{B_{23}}\ast E^{A_3}_{B_{23}})$ for $(ii)$. This follows from proposition \ref{exel} for $(i)$ and from proposition \ref{redfreeprod} for $(ii)$.  An immediate verification shows that the diagrams that need to commute do so, hence in both cases the generalized free product of the minimal triangle is realizable.
\rule{5 pt}{5 pt}

\par For unital  algebras, *-algebras and C*-algebras, we can prove a similar reduction theorem as the one for groups. 

\begin{theorem}
Let $(A_i)_{i\in I}$ be a family of unital algebras (*-algebras, C*-algebras), let for each $i\in I$ $B_{ij}\subset A_i,j\in I, j\ne i$ be unital subalgebras (*-algebras, C*-algebras) such that $B_{ij}\stackrel{\phi_{ij}}{\simeq}B_{ji}$, with $\phi_{ij}^{-1}=\phi_{ji},\forall i,j\in I,i\ne j$. Let $B_i:=\bigvee_{j\in I,j\ne i}B_{ij}\subseteq A_i$. The following are equivalent:\\
(i) the minimal amalgam of the family $(B_i)_{i\in I}$ is realizable.\\
(ii)) the generalized free product of $(A_i)_{i\in I}$ with amalgamated subalgebras $(B_{ij})_{i,j\in I,i\ne j}$ is realizable.
\end{theorem}

{\it Proof} The proof follows the lines of the reduction theorem for groups.  The implication $(ii)\Rightarrow(i)$ is trivial. For $(i)\Rightarrow(ii)$, note that  for each of the considered categories there exists a well defined free product with amalgamation. Denote by $B$ the minimal realizable amalgam. In each of the cases consider the algebra 
 \[ A_0=\ast_{B,i\in I}(A_i\ast_{B_i}B) .\]
 We denote with the same $A_i,B_{ij}$ the image of these algebras inside $A_0$. Note that in $A_0$ we have $A_i\cap A_j=B_{ij}=B_{ji}$.
 Then $A=\bigvee_{i\in I}A_i\subseteq A_0$ contains isomorphic copies of the factors $A_i$, hence the generalized free product of $(A_i)_{i\in I}$ is realizable. Even more is true. Since $B_i\subset A_i\subset A$ and $B=\ast_{i\in I}(B_i,(B_{ij})_{j\ne i})$, it follows that $B\subset A$. As in the case of groups, it is then easy to check that $A_0=A$. For example, in the case of C*-algebras, we have
 \[ \overline{(span\{ x_1x_2...x_n|n\ge1, x_k\in A_{i_k}\ast_{B_{i_k}}B,i_k\in I, i_1\ne i_2\ne ...\ne i_n\})}^{\| \|}=A_0. \]
 But $B$ is generated by all the $B_i$'s, so the products $x_1x_2....x_n$ all belong to $A$, which is the C*-subalgebra generated by $(A_i)_{i\in I}$. Therefore
  
  \[ A=\overline{(span\{ x_1x_2...x_n|n\ge1, x_k\in A_{i_k},i_k\in I, i_1\ne i_2\ne ...\ne i_n\})}^{\| \|}=\]
  \[=\overline{(span\{ x_1x_2...x_n|n\ge1, x_k\in A_{i_k}\ast_{B_{i_k}}B,i_k\in I, i_1\ne i_2\ne ...\ne i_n\})}^{\| \|}=A_0, \]
and we are done.
\rule{5 pt}{5 pt}


\section{\label{freeness}Reduced free products in the presence of the minimal amalgam}
\par
The next theorem shows that it is actually possible to have a "reduced" version of the generalized free product, if the minimal amalgam exists. We present a natural way to build a state that relates the generalized free product to the freeness concept introduced by Voiculescu.
\begin{theorem}
Let $(A_i)_{i\in I}$ be a family of unital C*-algebras,  for each $i\in I$ let $B_{ij}\subset A_i,j\in I, j\ne i$ be unital subalgebras such that $B_{ij}\stackrel{\varphi_{ij}}{\simeq}B_{ji}$, with $\varphi_{ij}^{-1}=\varphi_{ji},\forall i,j\in I,i\ne j$. Let $B_i:=\bigvee_{j\in I,j\ne i}B_{ij}\subset A_i$. Suppose that the minimal amalgam $B$ is realizable and that there are conditional expectations $\phi_i:A_i\to B_i$,$\psi_i:B\to B_i, i\in I$, each having faithful GNS representation. Then there exists an algebra $A$ and a conditional expectation $\phi:A\to B$ such that :\\
(i) for each $i\in I$, $A$ there exists an injective $\ast$-homomorphism $\sigma_i: A_i\to A$ and $A$ is generated by $\cup_{i\in I}\sigma_i(A_i)$;\\
(ii) $\phi|_{A_i}=\phi_i$ and the family $(\sigma_i(A_i))_{i\in I}$ is free with amalgamation over $B$ in $(A,\phi)$.
\end{theorem}

{\it Proof.  } Fix $i\in I$ and consider the reduced free product $(A_i\ast_{B_i}B,\phi_i\ast \psi_i)$. Since both $\phi_i$ and $\psi_i$ have faithful GNS representations, there exists a conditional expectation $E_i:A_i\ast_{B_i}B\to B$ like in prosposition \ref{bl}. Let $\pi_i:A_i\ast_{B_i}B\to \mathcal{L}(L^2(A_i\ast_{B_i}B,E_i))$ be the associated $\ast$-representation, as in section \ref{s0}. As described in \cite{DyBl},  $\pi_i|_{A_i}$ is faithful. Let
\[ (A,\phi )= \ast_{B,i\in I}(A_i\ast_{B_i}B,E_i) .\]
This is done via the free product representation $\lambda=\ast _i\lambda_i$ as described in section \ref{s0}. Because $\pi_i|_{A_i}$ is faithful and $\lambda_i=V_i(\pi_i\otimes Id)V_i^{*}$ where $V_i$ is the unitary from section \ref{s0}, we conclude that $\lambda_i|_{A_i}$ is also faithful.   Further more, $\lambda$ is precisely the GNS representation of $A$ associated to $\phi$. We have:
\[ A= \overline{\cup\lambda_i(A_i\ast_{B_i}B)}^{\|\cdot \|}=\]
\[ \overline{ span\{\lambda_{i_1}(x_1)\lambda_{i_2}(x_2)...\lambda_{i_n}(x_n)|x_j\in A_{i_j}\ast_{B_{i_j}}B,i_1\ne i_2\ne...i_n,n\ge 1\} }^{\|\cdot \|}=\]
\[ \overline{span\{\lambda_{i_1}(a_1)\lambda_{i_2}(a_2)...\lambda_{i_n}(a_n)|x_j\in A_{i_j},i_1\ne i_2\ne...i_n,n\ge 1\}}^{\|\cdot \|},\]
 because the copy of $B$ inside $A$ is generated by $\cup_{i\in I}\lambda_i(B_i)\subset \cup_{i\in I}\lambda_i(A_i)$.
 Letting $\sigma_i=\lambda_i|_{A_i}$ we get $(i)$. Remark that $(ii)$ is just a consequence of the construction; inside $A$, the algebra generated by $\lambda(A_i)$ and $\lambda(B)$ is exactly $\lambda(A_i\ast_{B_i}B)=\lambda_i(A_i\ast_{B_i}B)$.
 \rule{5 pt}{5 pt}
 \begin{remark}
We can actually describe $L^2(A,\phi)$ using the description of $L^2(A_i\ast_{B_i}B,E_i),i\in I$, given in \cite{DyBl}. For $i\in I$, let $H_{i,1}=L^2(A_i,\phi_i)$ and $H_{i,2}=L^2(B,\psi_i)$;  let $H_{i,k}^{\circ}$ denote the orthogonal complement with respect to $B_i$ (k=1,2). From \cite{DyBl} we get 
\[ F_i:=L^{2}(A_i\ast_{B_i}B,E_i)=B\oplus 
\bigoplus_{{\tiny \begin{array}{c}
 n\geq1\\
k_1,k_2,...,k_n\in \{1,2\}\\
k_1\not=k_2\not=...\not=k_n\ne 2
\end{array}}}
H^{\circ}_{i,k_1}\otimes_{B_i}H^{\circ}_{i,k_2}\otimes_{B_i}...\otimes_{B_i}H^{\circ}_{i,k_n}\otimes_{B_i}B.\]
By construction, we know the structure of $L^2(A,\phi)$ is (with $F^{\circ}_i$ the orthogonal complement with respect to $B$) 
\[ L^2(A,\phi)=B\oplus
\bigoplus_{ {\tiny \begin{array}{c}
 n\geq1\\
i_1,i_2,...,i_n\in I\\
i_1\not=i_2\not=...\not=i_n\ne 2
\end{array} }  }
F^{\circ}_{i_1}\otimes_{B}F^{\circ}_{i_2}\otimes_{B}...\otimes_{B}F^{\circ}_{i_n}=B\oplus\]

\[ 
\bigoplus_{{\tiny \begin{array}{c}
 n\geq1\\
i_1,i_2,...,i_n\in I\\
i_1\not=i_2\not=...\not=i_n\ne 2
\end{array}}}
\bigoplus_{{\tiny \begin{array}{c}
m_{i_1}\ge1,m_{i_2}\ge 1,...,m_{_n}\ge 1\\
k_{(i_p,1)},k_{(i_p,2)},...,k_{(i_p,m_p)}\in \{1,2\}, \\
k_{(i_p,1)}\ne k_{(i_p,2)}\ne...\ne k_{(i_p,m_p)}\ne 2,\\
p=\overline{1,n}
\end{array}}}
H^{\circ}_{i_1,k_{(i_1,1)}}\otimes_{B_{i_1}}H^{\circ}_{i_1,k_{(i_1,2)}}\otimes_{B_{i_1}}...\]
\[ \otimes_{B_{i_1}}H^{\circ}_{i_1,k_{(i_1,m_1)}}\otimes_{B_{i_1}}B\otimes_BH^{\circ}_{i_2,k_{(i_2,1)}}\otimes_{B_{i_2}}H^{\circ}_{i_2,k_{(i_2,2)}}\otimes_{B_{i_2}}...\]
\[H^{\circ}_{i_{n-1},k_{(i_{n-1},m_{n-1})}}\otimes_{B_{i_{n-1}}}B\otimes_BH^{\circ}_{i_n,k_{(i_n,1)}}\otimes_{B_{i_n}}...H^{\circ}_{i_n,k_{(i_n,m_n)}}\otimes_{B_{i_n}}B.\]

This description of the Fock space associated to $\phi$, and hence of the kernel of $\phi$, corresponds to the following intuitive decomposition of a word over the alphabet $\cup_{i\in I}A_i$ in the generalized free product with amalgamations. If $w=a_1a_2...a_n$, with $n\ge 1$ and $a_j\in A_{k_j},i_1\ne k_2\ne...\ne k_n$, then $a_j^{\circ}=a_j-\phi_{k_j}(a_j)\in A_{k_j}^{\circ}\subset H^{\circ}_{k_j,1}$, so between elements from the $A_i^s$ with $0$ expectation there will be elements from $B$, hence
\[ w\in B+\sum_{{\tiny \begin{array}{c}
m\ge 1\\
i_1, i_2,...i_m\in I
\end{array}}}A^{\circ}_{i_1}B(i_1,i_2)A^{\circ}_{i_2}B(i_2,i_3)...B(i_{m-1},i_m)A^{\circ}_{i_m},\]
where if $s\ne r$ then $B(r,s)$ is the set of words over $\cup_{i\in I}B_i$ that do not start with a "letter" in  $B_r$ and do not end with a "letter" from  $B_s$ and $B(r,r)=B^{\circ,r}\subset H^{\circ}_{r,2}$.

\end{remark}

\section{\label{examples}Examples of triangles of algebras}
\par
In this section we conclude our paper with 3 examples of realizable minimal amalgams, sharing some common features. One of this is that they all appear as the generalized free product of a family $(D_i,(D_{ij})_{j\ne i})_{i\in I}$, where $I=\{ 1,2,3\}$. The relations among the algebras can be schematically described using a triangle diagram of algebras and injective $\ast$-homomorphisms, represented by arrows:

\triunghi{D_1}{D_2}{D_3}{D_{12}}{D_{13}}{D_{23}}{D}

\begin{example}
 The first example presents the free amalgamated product of 3 unital C*-algebras as the result of the minimal amalgam of a triangle of unital C*-algebras. For this, let $A_1,A_2,A_2$ be unital C*-algebras and let $B$ be a common unital C*-subalgebra. Let 
 \[ A=\ast_{B,i}A_i .\]
 Let $D_i=A_j\ast_B A_k\subset A$, for $i,j,k\in \{ 1,2,3\}, i\ne j\ne k \ne i$. With $D_{ij}:=D_i\cap D_j=A_k$ ,$i,j,k\in \{ 1,2,3\}, i\ne j\ne k \ne i$ we get $D=B$ and it is an easy exercise to check that 
 \[ A=\ast_{i=1}^3(D_i,(D_{ij})_{j\ne i}). \]
 Suppose $\phi_i:A_i\to B$ are conditional expectations with faithful GNS representations. Then the triangle embeds also in the reduced amalgamated free product $(A,\phi)=\ast_i(A_i,\phi_i)$. By \ref{bl}, there exist conditional expectations which reverse the arrows of the triangle. Moreover, if we let $E_{ij}:D_i\to D_{ij}$ be this maps, then $E_{ij}|_{D_{ik}}=\phi_k$ for all $i1\le i\ne j\ne k\le2$.  We could think of $(A,\phi)$ as the reduced version of the minimal amalgam of the triangle.
\end{example}

\begin{example}
The second example presents the tensor product of 3 unital *-algebras ( C*-algebras ) as the minimal amalgam of a triangle. Let $A_1,A_2,A_3$ be unital $\ast$-algebras and let 
\[ A=A_1\otimes A_2\otimes A_3 .\]
Choose $D_i=A_j\otimes A_k$ , $D_{ij}=A_k=D_{ji}$ where $i,j,k\in \{ 1,2,3\}, i\ne j\ne k \ne i$. Automatically we need $D=\co$.  The injective homomorphisms from $D_{ij}=A_k$ into $D_i$ are the obvious ones, sending $A_k$ into its corresponding component in $D_i$. Note that the entire family embeds naturally into $A$, and $A$ is generated by the images of $D_i$, hence
\[ A=A_1\otimes A_2\otimes A_3=\ast_{i=1}^3(D_i, (D_{ij})_{j\ne i}).\]
Suppose now that $\phi_i:A_i\to \co$ are states. We can then define conditional expectations $E_{ij}:_i\to D_{ij}$, letting $E_{ij}(a_k\otimes a_j)=\phi_k(a_k)a_j$. Then $E_{ij}|_{D_{ik}}=\phi_k$.
\end{example}

\begin{example}
Denote by $M_n:=M_n(\co)$ the algebra of n by n matrices of complex numbers. The previous example shows the trivial way to get $M_8\simeq M_2\otimes M_2\otimes M_2$ as the minimal amalgam of a triangle. In this example we show that there is another way of getting $M_8$ as the minimal amalgam of a triangle with the same "vertices" as above. Let $D_i=M_2\otimes M_2$ for $i=1,2,3$. Consider the following two selfadjoint unitaries :
\[ u=\left( \begin{array}{cccc}
1&0&0&0\\
0&0&0&1\\
0&0&1&0\\
0&1&0&0
\end{array}\right),
v=\left( \begin{array}{cccc}
1&0&0&0\\
0&1&0&0\\
0&0&0&1\\
0&0&1&0
\end{array}\right) \]

It is worth mentioning that $u$ and $v$ are the only biunitaries in $M_4$ which are also permutation matrices; a matrix $w=(w^a_b)^i_j\in M_n\otimes M_k$ is called a biunitary if  it is a unitary matrix and if the matrix defined by $((w_1)^a_b)^i_j=(w^b_a)^i_j$ - the block-transpose of $w$ - is also a unitary (see \cite{sun} for more details on biunitary permutation matrices).
Since $u,v$ are biunitaries, the following two are non-degenerated commuting squares :  \\
\begin{table}[htbp]
\begin{center}
\begin{tabular}{ccccccc}

$u(M_2\otimes I_2)u$&$\subset$&$M_2\otimes M_2$&$ $&$v(M_2\otimes I_2)v$&$\subset$&$M_2\otimes M_2$\\
$\cup$&$ $&$\cup$&$   ,   $&  $\cup$&$ $&$\cup$\\
$\co$&$\subset$&$I_2\otimes M_2$&$ $&$\co$&$\subset$&$I_2\otimes M_2$

\end{tabular}
\end{center}
\end{table}
 The non-degeneracy means that, if we let $ D_{12}=u(M_2\otimes I_2)u, D_{23}=v(M_2\otimes I_2)v, D_{13}=I_2\otimes M_2$, then 
 \[ M_2\otimes M_2=span\left\{ AB | A\in D_{12},B\in D_{13}\right \}=span\left\{ BA | A\in D_{12},B\in D_{13}\right \},\]
 and also  
\[ M_2\otimes M_2=span\{ AB | A\in D_{23},B\in D_{13}\}=span\{ BA | A\in D_{23},B\in D_{13} \}.\]
 Even more is true in commuting squares : let $\tau_n$ to be the normalized trace in $M_n$ and define $E_0:M_2\otimes M_2\to I_2\otimes M_2$,  $E_0(a\otimes b)=\tau_2(a)(I_2\otimes b)$;  $E_u:M_2\otimes M_2\to D_{12}$ , $E_u(a\otimes b)=\tau_2(b)u(a\otimes I_2)u$ and finally $E_v:M_2\otimes M_2\to D_{23}$ , $E_v(a\otimes b)=\tau_2(b)v(a\otimes I_2)v$. Then $E_u(I_2\otimes b)=\tau_2(b)=E_v(I_2\otimes b)$. Let $e(i,j)$ be the matrix units of $M_2$ and let $e_u(i,j)=u(e(i,j)\otimes I_2)u$ , $e_v(i,j)=v(e(i,j)\otimes I_2 )v$ and $e_0(i,j)=I_2\otimes e(i,j)$.  One checks that we also have  
 \[ M_2\otimes M_2=span\{ AB | A\in D_{12}, B\in D_{23}\}=span\{ BA | A\in D_{12}, B\in D_{23}\},\] and that 
 $D_{12}\cap D_{23}=\co$.
However, $E_u(v(a\otimes I_2)v)=a_{11}e_u(1,1)+a_{22}e_u(2,2)$, so that the commuting squares condition fails for this pair of algebras. 
We want to find the minimal amalgam. For this purpose note the relations :
\begin{table}[htdp]
\caption{Commutation relations}
\begin{center}
\begin{tabular}{lcl}
$e_u(i,i)x_v=x_ve_u(i,i),$&$\mbox{ for all }i\in I,x_v\in D_{23}$\\
$e_u(i,i)e_0(k,k)=e_0(k,k)e_u(i,i)$&$\mbox{ for all }i,k$&$$\\
$e_u(i,j)e_0(k,l)=e_0(k,l)e_u(\sigma(i),\sigma(j)),$&$ \mbox{ for all } i,j,k,l$&$$,\\

$e_v(i,i)x_0=x_0e_v(i,i),$&$\mbox{ for all }i\in I,x_0\in D_{13}$\\
$e_v(i,i)e_u(k,k)=e_u(k,k)e_v(i,i)$&$\mbox{ for all }i,k$&$$\\
$e_v(i,j)e_u(k,l)=e_u(k,l)e_v(\sigma(i),\sigma(j)),$&$ \mbox{ for all } i,j,k,l$&$$,\\

$e_0(i,i)x_u=x_ue_0(i,i),$&$\mbox{ for all }i\in I,x_v\in D_{12}$\\
$e_0(i,i)e_v(k,k)=e_v(k,k)e_0(i,i)$&$\mbox{ for all }i,k$&$$\\
$e_0(i,j)e_v(k,l)=e_v(k,l)e_0(\sigma(i),\sigma(j)),$&$ \mbox{ for all } i,j,k,l$&$$,\\

\end{tabular}
\end{center}
\label{co}
\end{table}

 $\sigma=(12)\in S_2$. As a vector space, because of this commutation relations (see \ref{dmodul}), the minimal amalgam $A$ will be generated by the non-zero elements 
 \[ \{e_u(i,j)e_v(k,l)e_0(r,s)| i,j,k,l,r,s\in I\}.\]  
 
 
 The products of minimal projections generate in the amalgam 8 minimal projections, $\{ e_u(i,i)e_v(k,k)e_0(r,r) | i,k,r=1,2 \}$, which add up to one; the partial isometries can then be constructed and a complete system of matrix units for $A$ can be found. Even more, the "vertices" $M_2\otimes M_2$ embed in this algebra in such a way that the resulting diagrams commute. Let us check it for $D_2$. Define $\phi_2:D_2\to A$ by

 \[ \phi(a\otimes b)=\phi_2\left( \sum_{i,j,k,l=1}^{2}x(i,j,k,l)e_u(i,j)e_v(k,l)\right)=\sum_{i,j,k,l=1}^{2}x(i,j,k,l)e_u(i,j)\cdot e_v(k,l)\cdot1.\]
 It is a straitforward computation that shows the one-to-one correspondence between $\{ e_u(i,j)e_v(k,l) |1\le i,j,k,l\le 2 \}$ and $\{ e_4(\alpha,\beta) | 1\le \alpha,\beta\le 4\}$- the matrix units of $M_4$. This shows that $\phi_2$ is a well defined linear map. It is easy to check that it preserves conjugation. We will proceed to prove injectivity. Suppose $\phi_2(\sum_{i,j,k,l=1}^{2}x(i,j,k,l)e_u(i,j)e_v(k,l))=0$. Then 
 \begin{equation}\label{suma}
  \sum_{i,j,k,l=1}^{2}x(i,j,k,l)e_u(i,j)\cdot e_v(k,l)\cdot 1=0 \mbox{ in } A.
  \end{equation}
 Suppose there are $i_0,j_0,k_0,l_0$ such that $x(i_0,j_0,k_0,l_0)\ne 0$. Multiply \ref{suma} by the left with $e_u(i_0,i)$ and by the right with $e_u(j,j_0)$ to get
 \[ \sum_{i,j,k,l=1}^{2}x(i,j,k,l)(e_u(i_0,i)e_u(i,j)e_u(j,j_0))\cdot e_v(\sigma (k),\sigma(l))\cdot 1=0 \]
 or
 \begin{equation}\label{su2}
  \sum_{k,l=1}^{2}x(i_0,j_0,k,l)e_u(i_0,j_0)\cdot e_v(\sigma (k),\sigma(l))\cdot 1=0. 
\end{equation}

 Multiply again from the left with $e_u(j_0,i_0)$ to get
\[ \sum_{k,l=1}^{2}x(i_0,j_0,k,l)e_u(j_0,j_0)\cdot e_v(\sigma (k),\sigma(l))\cdot 1=0. \]
hence
\[ \sum_{k,l=1}^{2}x(i_0,j_0,k,l)e_u(i_0,j_0)\cdot e_v(k,l)\cdot 1=0. \]
Still in \ref{su2} multiply from the right with $e_u(j_0,i_0)$, use the fact that $\sigma^2=id$ to get
\[ \sum_{k,l=1}^{2}x(i_0,j_0,k,l)e_u(i_0,i_0)\cdot e_v(k,l)\cdot 1=0. \]
Multiply in both \ref{su2} and the last equation from the left with $e_v(\sigma(k_0),\sigma(k))$ and from the right with $e_v(l,l_0)$ and we get:
\[ e_u(i_0,i_0)\cdot e_v(k_0,k_0)\cdot 1=0, \]
\[ e_u(j_0,j_0)\cdot e_v(k_0,k_0)\cdot 1=0. \]
and similarly for $l_0$. Even if $i_0=j_0$ or $k_0=l_0$ we will still get $e_u(i,i)e_v(k,k)=0$ for $1\le i,k\le2$. Adding them up we get $1=0$ in the amalgam, a contradiction, since the amalgam contains at least $\co$. So $\varphi_2$ is injective. One can define similar maps $\phi_1:D_1\to A,\phi_3:D_3\to A$, which turn out to be injective. It is also clear that $\phi_1|_{D_{12}}=\phi_2|_{D_{12}}$,$\phi_1|_{D_{13}}=\phi_3|_{D_{13}}$, and $\phi_2|_{D_{23}}=\phi_3|_{D_{23}}$. Hence $A\simeq M_8$ is the realizable minimal amalgam of the triangle of algebras.
 \end{example}

\begin{remark}
One can generalize the previous example to arbitrary unitaries and vertex algebras of arbitrary dimensions. With the hypothesis that the commuting squares of algebras appearing in the  construction are non-degenerated (i.e. the vertex algebras are the span of the product of the edge algebras) and that the edge algebras intersect along the complex numbers, one can see that the result of the amalgam will always be the span of the product of the edge algebras.
\end{remark}

\begin{remark}The reduction theorems presented in the paper point to the fact that the minimal amalgam plays the key role. We are still looking for a good sufficient condition for such a minimal amalgam to exist. 
\end{remark}
\begin{acknowledgements} I would like to thank for the kind hospitality received at I.H.P. in Paris during the special semester on K-theory and non-commutative geometry, when part of this work was done. Special thanks go to prof. Ryszard Nest and Etienne Blanchard for discussions, suggestions and encouragements. I want also to thank my advisor, prof. Florin Radulescu, for the constant help and guidance. 
\end{acknowledgements}

 \end{document}